# Iterative Designs with Similar Triangles

Chris Thron, Texas A&M University-Central Texas

One of the best things about geometry is that it's *cool*! Geometry enables us to create incredible designs and astounding patterns.

This article shows how to use a simple technique (iteration) to create designs that are both cool *and* introduce important mathematical concepts. We give examples using one of the most common of geometrical constructions, namely, the similar triangle. In the final section, we suggest other possibilities that could be used for guided exploration or for student projects.

## Iteration, Dynamic Geometry Software, And Macros

One powerful method for creating intricate yet symmetrical designs is *iteration*, that is, applying the same geometrical construction repeatedly in ways that give rise to additional patterns. Modern interactive geometry software is perfect for constructing fantastic iterative designs. The pictures in this article were created using "Compass and Ruler" ([car.rene-grothmann.de/](car.rene-grothmann.de/) ), a friendly and versatile Web-based Java program authored and maintained by René Grothmann. These constructions can be reproduced in other geometry programs: a list of alternatives may be found at [en.wikipedia.org/wiki/List_of_interactive_geometry_software](en.wikipedia.org/wiki/List_of_interactive_geometry_software). One particularly popular alternative is GeoGebra ([www.geogebra.org/](www.geogebra.org/)).

The reason that dynamic geometry software is so well suited for creating iterative designs is that repetitive procedures can be automated through the creation of *macros*. Macros function like small computer programs, although they do not have to be programmed in a computer language. Instead, the user performs a series of actions on a set of points and lines, and instructs the computer to record these actions (which may include coloring various picture elements, or hiding accessory elements that are not wanted in the final picture). The computer is then able to repeat the same steps at will using a different set of points and lines. (Readers familiar with the "macro" features in Microsoft Excel or Word will note the similarity).

A simple example of a macro is shown in Figure 1. Starting with three points shown at left, the user went through the steps of joining the points to create a triangle; finding the midpoints of the triangle's sides; and joining the midpoints to create an inner triangle. The user then let the computer take a "snapshot" of the final construction, saving it as the macro "MidTri."

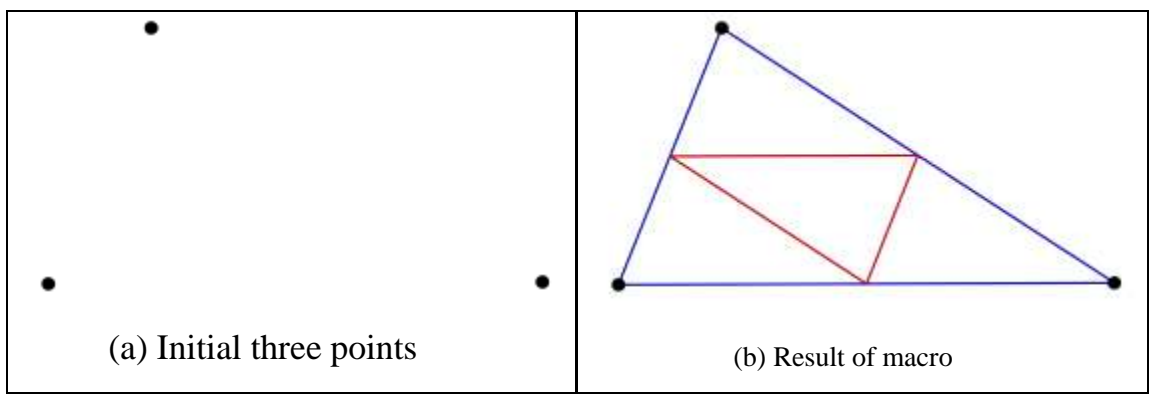

(a) Initial three points

(b) Result of macro

**Figure 1** The macro "MidTri" applied to three initial points *(a)* creates a triangle joining the midpoints *(b)*.

In C.a.R. (as in other dynamic geometry programs) it is possible to move the original points around, while preserving the construction. A student who plays with the three original points may notice the striking fact that the inner triangle always seems to be similar to the outer triangle, and exactly half the size. This of course may be proved by the usual methods.

Once the macro has been created, the user can then invoke it at any time and apply it to any three points at will. The design shown in Figure 2 was created by successively applying the macro to each new triangle created by the previous application of the macro.

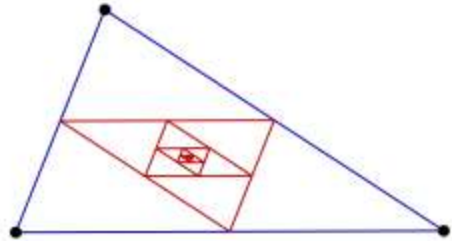

**Figure 2** The macro in Figure 1 can be iterated to create a decreasing series of triangles.

## Iterated triangles and concurrencies

The iteration in Figure 2 can lead to an important geometrical discovery. To bring this out, we create another macro "RGBMidTri" that builds on MidTri by adding the red, green, and blue median segments shown in Figure 3.

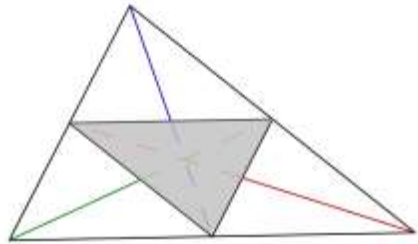

**Figure 3** Macro RGBMidTri adds median segments to the midpoint triangle created by MidTri.

Iterating RGBMidTri produces two green, two red, and two blue segments that approach closer and closer to a common point, as shown in .Figure 4.

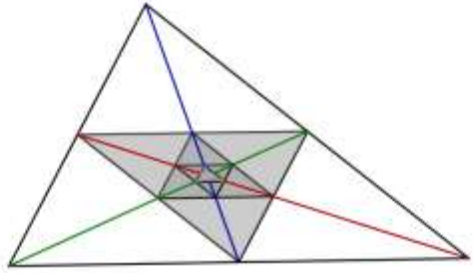

**Figure 4** When RGBMidTri is iterated, all the median segments approach the centroid.

Using the fact that each triangle in the series is exactly half the size of the previous one, by comparing the lengths of green segments added in alternate iterations one may verify that the two green segments' lengths have ratio 2:1; and similarly for the blue and red segments. This illustrates the Median Concurrency Theorem, namely:

> The three medians of a triangle meet at a single point (*called the centroid*). This point divides each of the medians into two segments whose lengths have the ratio 2:1.

By adding altitudes and perpendicular bisectors to this iterated diagram, it is possible to produce a construction that visually demonstrates the *Euler Line Property*:

> The centroid of a triangle lies on the segment that connects the triangle's circumcenter (*intersection of perpendicular bisectors*) and orthocenter (*intersection of altitudes*), and divides the segment in the ratio 2:1.

For a pictoral representation, see Figure 11 at the end of this article.

### Inscribing and Circumscribing Triangles

The above iteration suggests other possibilities. The similarity of the inscribed triangles naturally suggests the question: are there other ways to inscribe a similar triangle inside a similar triangle? A little experimentation with cutouts of similar triangles of different sizes provides convincing evidence that there are other ways as well, as shown in Figure 5. But although this may be easy to do with cut-outs, it's a bit more difficult on the computer.

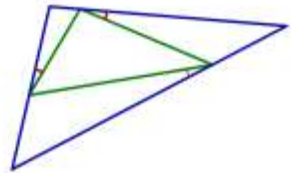

**Figure 5** A green similar triangle is inscribed in the blue triangle by rotating at a particular angle.

The two triangles in Figure 5 are rotated relative to each other, where each angle marked in red in the figure is equal to the angle of rotation. This picture suggests a way to *circumscribe* a similar triangle – namely, by drawing lines at equal angles from the three vertices of the starting triangle. A macro can be created to do this, as shown in Figure 6. The macro RGBCircumTri takes a given angle (determined by three points) and a given triangle to create a circumscribing triangle with red, green, and blue sides that has relative rotation equal to the given angle.

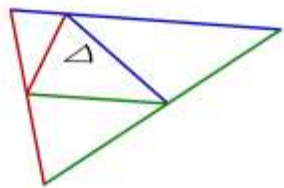

**Figure 6** The C.a.R. macro RGBCircumTri circumscribes a similar triangle that has a given relative angle of rotation (the black angle inside the inner triangle).

**Spirals, Limits, and Perspective**

Now the real fun starts! Using the macro shown in Figure 6, we can iteratively circumscribe triangles around triangles. The triangles become larger and larger, so we need to zoom out as we create more triangles (C.a.R. provides this capability). Two possibilities are shown in Figure 7: the different figures correspond to different angles of rotation between successive triangles.

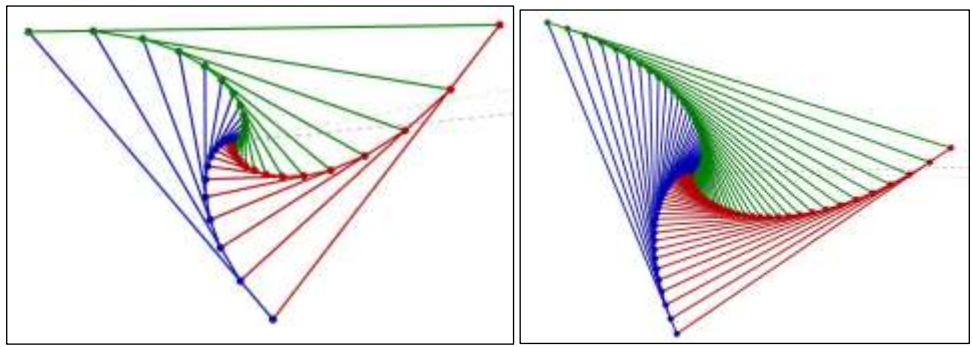

**Figure 7** Iteration of the macro RGBCircumTri produces similar spirals for two different rotation angles.

Instead of triangles, we see spirals! Each successive circumscribe operation creates a new triangle that always has the same size ratio and the same relative angle with the previous triangle. Thus the size of the triangles increases in proportion to the angle of rotation – producing *logarithmic* spirals.

If we zoom back in and change the shape of the first triangle by moving the vertices around, the entire pattern will change as a result. But no matter what changes we make, the three spirals will always have similar shapes, even if the triangle has sides that are quite unequal. The more unequal the sides, the "flatter" the spirals will be.

When zoomed all the way out, the inner triangles appear to shrink to a point – this illustrates the mathematical concept of limit. Students may discover the surprising fact that the limit points remains the same even if the angle of rotation between triangles is changed. This point actually has a name: it's called a "Brocard point," named after Henri Brocard who investigated them in 1875. We will return to this point later on in the article.

Zooming in or out (holding down the "+" or "-" buttons in C.a.R.) gives the illusion of traveling into or out of a rotating tunnel. This is because the point of convergence of the spirals also serves as a vanishing point for one-point perspective. It appears that the tunnel is rotating at a constant rate, because a triangle's angle of rotation is directly proportional to its size.

## Circles, Angles, and Equal Intercepted Arcs

What happens if we circumscribe several different similar triangles around the same triangle? Triangles of different colors are circumscribed around the red triangle in Figure 8. Strangely, it seems that corresponding vertices seem to lie on a circle. We can verify this by drawing a circle through three of the vertices (alternatively, C.a.R. has a tool for drawing a conic section through any five given points).

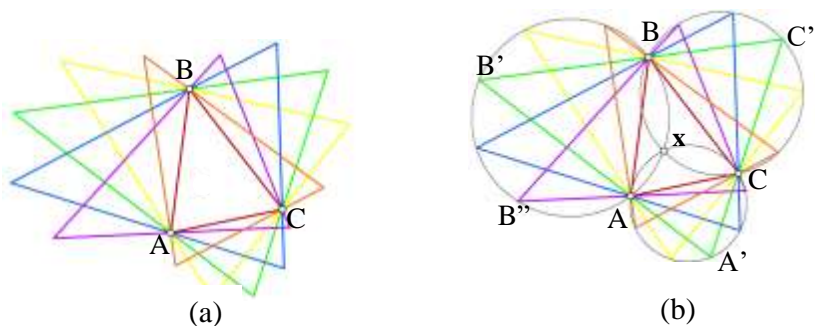

(a)          (b)

**Figure 8** In (a), the colored triangles are circumscribed around the red triangle at different angles of rotation. In (b), circles are drawn through the corresponding vertices of the circumscribed triangles.

There is a reason behind the circles! In Figure 8(b), notice for instance that green ∠B'

and purple ∠B" are both congruent to angle ∠B of the original (red) triangle: furthermore, both angles intercept arc AB on the leftmost circle. A well-known property of circles states that

> Any two angles inscribed in a circle that intercept the same arc are congruent.

The *converse* is also true, namely

> Any two congruent angles that intercept the same segment are both inscribed in a circle that also contains the endpoints of the segment.

This statement can be seen from Figure 9. If an angle's vertex lies inside or outside the circle determined by the other vertex and the endpoints of the segment, then it subtends a larger or smaller angle.

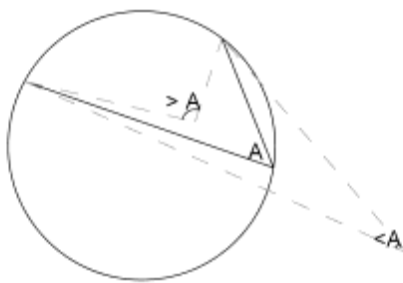

**Figure 9** Vertex points inside and outside the circle create angles that are bigger or smaller than vertex points on the circle, if all angles intercept the same arc.

There is a remarkable "coincidence" in Figure 8(b): the three circles apparently intersect at a single point **x**. By moving the points A, B, and C, the student may discover that this is always true, no matter what the shape of the original triangle. The reason may once again be found in the properties of angles inscribed in circles. If x lies on circle AB'B, it follows that the angle ∠AxB must be supplementary to ∠AB'B, which is equal to ∠B of the original triangle. By the same reasoning, if x also lies on circle BC'C then ∠BxC must be supplementary to ∠C. This forces ∠AxC to be supplementary to ∠A, which in turn forces **x** to lie on circle CA'A.

One final remarkable property of x is that it forms equal angles with all sides of the original triangle, so that ∠xAC, ∠xCB, and ∠xBA are all congruent. Furthermore, x has the same property relative to *any* of the colored triangles in Figure 8(b). The proof, once again, depends on successive applications of the fact that angles inscribed in the same circle that intercept the same arc are congruent. For instance, ∠xAC and ∠xBA intercept the same arc on the leftmost circle in Figure 8(b).

The iterative circumscribing of similar triangles in Figure 7 can be combined with the circles in Figure 8(b) to create the figure shown in Figure 10. There are lots of similarities in this figure! The figure shows that the Brocard point that is the point of convergence in Figure 7 is actually one and the same as the three-circle intersection shown in Figure 8(b).

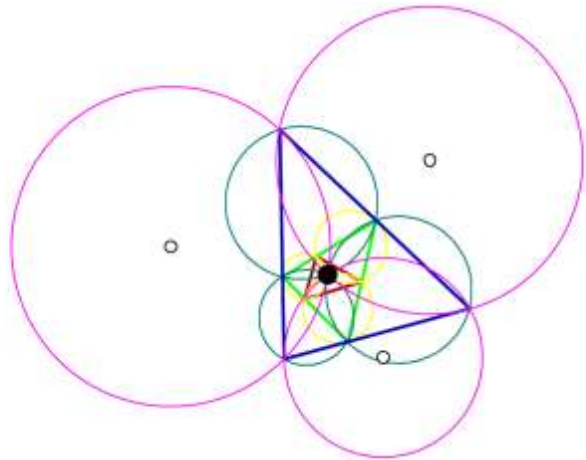

**Figure 10** All circles intersect at the same point!

## Further Explorations

In this article we have shown just some of the design possibilities for iteration of a simple geometrical construction. There are many more that creative students (and teachers!) can find out for themselves. Here we offer a few suggestions for further explorations:

- The macro RGBCircumTri produces a circumscribing triangle whose sides make equal *clockwise* angles from the inner triangle's three sides. Suppose we measure equal *counterclockwise* angles instead?

- There's yet another family of circumscribed similar triangles around a given triangle. Can you discover it, and discover its properties? (Hint: one member of the family is represented in Figure 1).

- The iterative circumscribing process shown in Figure 7 can be performed on other polygons besides triangles. Are the circumscribing polygons still similar to the original? Is there still a unique point of convergence?

- Other geometrical elements besides triangles can be used as the basis of iteration. Circles can also be used to create some interesting optical-illusion effects.

## Web Resources

Bach, Michael, "87 Optical Illusions & Visual Phenomena," http://www.michaelbach.de/ot/ .

Grothmann, René "C.a.R. Documentation,"

  http://zirkel.sourceforge.net/doc_en/Documentation/index.html

Weisstein, Eric W. "Brocard Points." From MathWorld--A Wolfram Web Resource.

  http://mathworld.wolfram.com/BrocardPoints.html

## BIBLIOGRAPHY

Serra, Michael *Discovering Geometry: an Investigative Approach, Fourth Edition*. Emeryville,

CA: Key Curriculum Press, 2008.

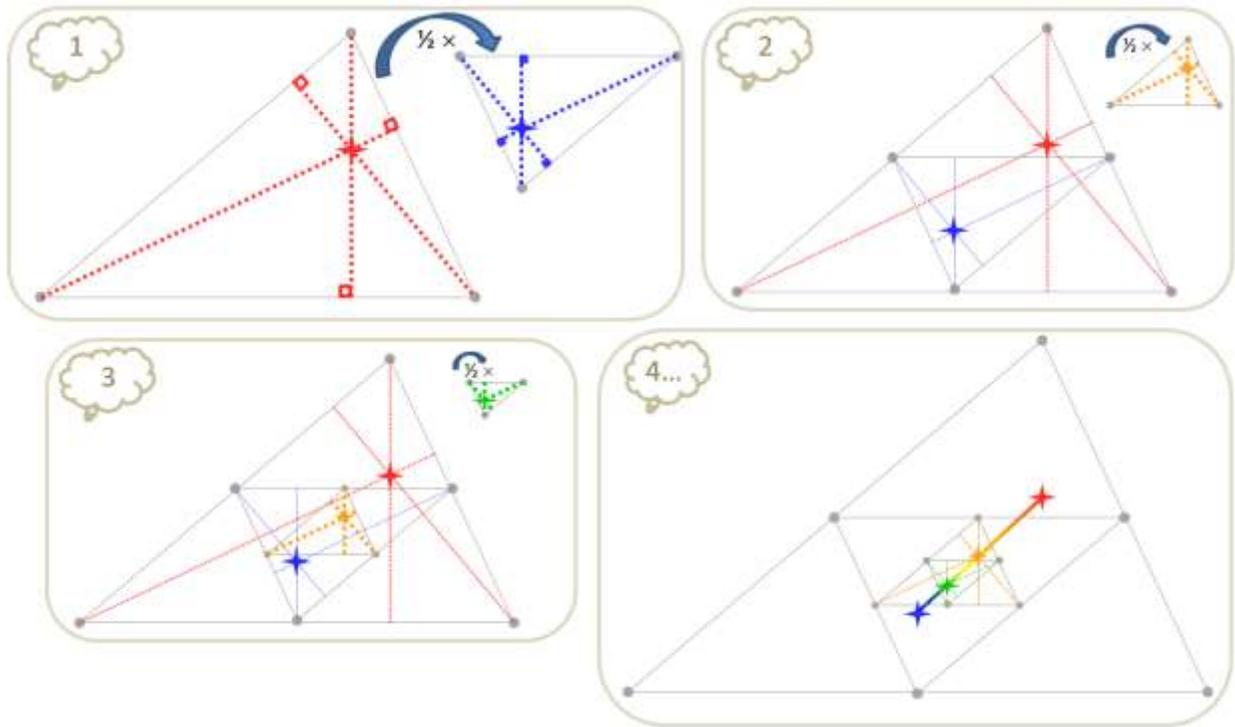

Figure 11 Illustration of the Euler line, using iteration. The orthocenter (altitude concurrency point) of each triangle is the circumcenter (perpendicular bisector concurrency point) of the next larger triangle. By parallelism, we know that all points lie on the same line. Since the triangles converge down to the median, the median also lies on the same line. All similar lengths in successive triangles are in the ratio of 2:1.